\documentclass[seceqn,dvips]{arxbj}

\aid{0}
\volume{16}
\issue{2}
\pubyear{2010}
\firstpage{418}
\lastpage{434}
\doi{10.3150/09-BEJ220}

\makeatletter
\newtheorem{theorem}{Theorem}

\makeatother

\begin{document}
\begin{frontmatter}

\title{Strong approximations of level exceedences related to multiple
hypothesis testing}
\runtitle{Strong approximations of exceedences}

\begin{aug}
\author[1]{\fnms{Peter} \snm{Hall}\thanksref{1}\ead[label=e1]{halpstat@ms.unimelb.edu.au}}
\and
\author[2]{\fnms{Qiying} \snm{Wang}\thanksref{2}\ead[label=e2]{qiying@maths.usyd.edu.au}\corref{}}
\runauthor{P. Hall and Q. Wang}
\address[1]{Department of Mathematics and Statistics, The University
of Melbourne, VIC 3010, Australia.\\ \printead{e1}}
\address[2]{School of Mathematics and Statistics F07, University of
Sydney, NSW 2006, Australia.\\ \printead{e2}}
\end{aug}

\received{\smonth{10} \syear{2008}}
\revised{\smonth{6} \syear{2009}}

%
\begin{abstract}
Particularly in genomics, but also in other fields, it has become commonplace
to undertake highly multiple Student's $t$-tests based on relatively
small sample sizes.
The literature on this topic is continually expanding, but the main
approaches used to control
the family-wise error rate and false discovery rate are still based on
the assumption that
the tests are independent. The independence condition is known to be
false at the level of
the joint distributions of the test statistics, but that does not
necessarily mean, for the
small significance levels involved in highly multiple hypothesis
testing, that the assumption
leads to major errors. In this paper, we give conditions under which
the assumption of
independence is valid. Specifically, we derive a strong approximation
that closely links
the level exceedences of a dependent ``studentized process'' to those
of a process of
independent random variables. Via this connection, it can be seen that
in high-dimensional,
low sample-size cases, provided the sample size diverges faster than
the logarithm of the
number of tests, the assumption of independent $t$-tests is often justified.
\end{abstract}

%
\begin{keyword}
\kwd{false discovery rate}
\kwd{family-wise error rate}
\kwd{genomic data}
\kwd{large deviation probability}
\kwd{moving average}
\kwd{Poisson approximation}
\kwd{Student's $t$-statistic}
\kwd{upper tail dependence}
\kwd{upper tail independence}
\end{keyword}

\end{frontmatter}

\section{Introduction}\label{sec1}

 Today it is commonplace to undertake highly multiple hypothesis
testing, generally in genomics and very often using tests based on
Student's $t$-statistic; see, for example, Benjamini and Yekutieli
(\citeyear{2001Benjamini}), Efron and Tibshirani (\citeyear{2002Efron}), Cui and Churchill (\citeyear{2003Cui}),
Amaratunga and Cabrera (\citeyear{2004Amaratunga}), page~114, Scheid and Spang (\citeyear{2005Scheid}),
Shaffer (\citeyear{2005Shaffer}), Fox and Dimmic (\citeyear{2006Fox}), Hu and Willsky (\citeyear{2006Hu}), Qiu and
Yakovlev (\citeyear{2006Qiu}), Efron (\citeyear{2007aEfron}), Liu and Hwang (\citeyear{2007Liu}) and van de Wiel
and Kim (\citeyear{2007Wiel}). This popularity of multiple $t$-testing also extends to
other fields (e.g., Pawluk-Ko\l c \textit{et al.}~(\citeyear{2006Pawluk})). The
principal methods used to control the family-wise error rate and false
discovery rate are founded on the assumption of independence among
tests. Alternative approaches are generally based either on Bonferroni
bounds, which are unsatisfactory for a variety of reasons (see, e.g.,
Perneger (\citeyear{1998Perneger})), or on the hope that, despite ample evidence
of non-independence in terms of correlation analysis, independence can
be assumed in practice.

The latter hope tends to be pinned either on work of Benjamini and
Yekutieli (\citeyear{2001Benjamini}), who argued that in some settings, the absence of
independence can give conservative results, or on experience with the
analysis of financial data, which suggests that in some circumstances,
it might be reasonable to assume that the upper tails of the test
statistics are independent, even if the joint distributions are not.
Upper tail independence, as it is sometimes called (for discussion,
see, e.g., Wu (\citeyear{1994Wu}), Falk and Reiss (\citeyear{2001Falk}), R.~Schmidt (\citeyear{2002Schmidt}),
Li (\citeyear{2006Li}), R.~Schmidt and Stadtm\"{u}ller (\citeyear{2006Schmidt}), T.~Schmidt (\citeyear{2007Schmidt})),
is generally assumed to be non-asymptotic in nature. That is, tails of
joint distributions are often taken to be perfectly independent beyond
a certain threshold.

However, this type of model is not really appropriate for the analysis
of genomic data. In particular, it is difficult to determine a
biological reason for, or the actual location of, a threshold. It is of
greater practical interest to consider the possibility that the
strength of dependence in upper tails could become successively weaker
as the number of simultaneous tests, and the number of data vectors,
increases. If this could be established in the context of tests based
on Student's $t$-statistic, it would lend immediate justification to
the often-made assumption (see the articles cited in the first
paragraph of this paper) that highly multiple $t$-statistics can be
taken to be independent.

The present paper will establish such a result. The mechanism for our
model involves the critical points for tests becoming more extreme as
the number, $p$, of tests diverges (in fact, the increase in critical
points is a direct consequence of $p$ diverging) so that the tests are
conducted further into the tails; furthermore, the tails of the
distributions of test statistics becoming successively lighter as the
number of degrees of freedom of the test statistics increases.

We impose particularly weak conditions on the marginal distributions of
components. In particular, the distributions need only three finite
moments. With this assumption, and permitting the size, $n$, of the
group sample to increase a little faster than the logarithm of the
number of tests, it follows from our results that the joint
distributions of test statistics enjoy an asymptotic form of the upper
tail independence property.

This result would not be so striking if the statistics had normal
distributions, but it fails for heavy-tailed distributions such as
those for which not all moments are finite. Of course, Student's
$t$-distribution is itself in this category, yet our results show that
asymptotic independence holds in a particularly strong sense for
Student's $t$-statistic, even if it is computed from relatively
heavy-tailed data. The reason this is possible is that we permit the
group sample size to increase at a rate that is just sufficient to
convert heavy tails to tails that are sufficiently light, to enable
approximate independence at high levels.

It can be seen from this property that the availability of upper-tail
asymptotic independence is a bonus of working with highly multiple
hypothesis testing, that is,~with ``large $p$ and small $n$'' problems.
It is not available in more conventional, ``small $p$ and large $n$''
problems, where there is a very large literature on modelling
dependence in highly multiple hypothesis testing.

There is a literature on comparing studentized means when the variances
used for studentizing are computed from pooled data and so are common
to each test statistic. However, in our experience, that approach is
used less frequently, in practice, than the ``local'' standardization
treated in the present paper. When using the latter method, each mean
is divided by the standard deviation of the sample from which it was
computed. A major motivation is that the true variances may be
different in each instance. Even if the variances can reasonably be
assumed to be the same, it can be desirable to use the local approach
since it confers greater robustness. For example, when applied to the
mean alone, rather than its locally studentized form, the
large-deviation properties that underpin the analysis of high-level
exceedences require the data to have lighter tails.

Statistical literature on highly multiple hypothesis testing is
outlined in helpful reviews by Hochberg and Tamhane (\citeyear{1987Hochberg}), Pigeot
(\citeyear{2000Pigeot}), Dudoit \textit{et al.}~(\citeyear{2003Dudoit}), Bernhard \textit{et al.}~(\citeyear{2004Bernhard})
and Lehmann and Romano (\citeyear{2005LehmannR}), Chapter 9. Benjamini and Hochberg (\citeyear{1995Benjamini})
introduced an approach, which has become very popular, to the
controlling of false discovery rates; see also Simes (\citeyear{1986Simes}), Hommel
(\citeyear{1988Hommel}), Hochberg (\citeyear{1988Hochberg}), Sarkar and Chang (\citeyear{1997Sarkar}),
Sarkar (\citeyear{1998Sarkar}), Sen
(\citeyear{1999Sen}), Hochberg and Benjamini (\citeyear{1990Hochberg}) and Lehmann \textit{et
al.}~(\citeyear{2005Lehmann}). Benjamini and Yekutieli (\citeyear{2001Benjamini}) specified conditions under
which simultaneous, dependent hypothesis tests, conducted as though
they were independent, give conservative results; Benjamini and
Yekutieli (\citeyear{2005Benjamini}) addressed similar issues in the context of false
coverage-statement rate. Sarkar (\citeyear{2002Sarkar}) extended the work of Benjamini
and Yekutieli~(\citeyear{2001Benjamini}). Efron (\citeyear{2007bEfron}) suggested correlation corrections
for large-scale simultaneous hypothesis testing. Blair \textit{et
al.}~(\citeyear{1996Blair}) proposed methods for controlling family-wise error rates in
multiple procedures, Holland and Cheung (\citeyear{2002Holland}) discussed robustness of
family-wise error rates and Clarke and Hall (\citeyear{2007Clarke}) discussed robustness
of testing procedures based on means.

\section{Results and applications}\label{sec2}

 \subsection{Model and main results}\label{sec2.1}

Given $p,n\geq1$, assume that for $1\leq i\leq p$ and $1\leq j\leq n$,
we observe data $U_{ij}$, which we\vspace{1pt} use to construct $t$-statistics $T_i=
n^{1/2} \bar{U}_i/S_i$, where $\bar{U}_i= n^{-1} \sum_j U_{ij}$
and $S_i^2=
n^{-1}
\sum_j U_{ij}^2-\bar{U}_i^2$.\vspace{1pt} In practice, the statistic $T_i$ is
used to test
the hypothesis that the $i$th group has zero mean, against a one-sided
alternative. When controlling the level of family-wise error rate
(FWER) for step-down tests, we require the values of probabilities
$P(T_i>t$ for $i=i_1,\ldots,i_k)$ for different levels $t$ and
different subsets $\{i_1,\ldots,i_k\}$ of $\{1,\ldots,p\}$.
Theorem~\ref{thm1}
below will enable us to compute these through approximation by the case
where the $T_i$'s are all independent; see Section~\ref{sec2.3} for further details.

We standardize $S_i^2$ by dividing by $n$, rather than $n-1$, since the
former is more common in nonparametric problems, but the results below
are unaffected by this issue. Since we studentize, there is no loss of
generality in assuming that the variance of each component equals~1.
More particularly, we ask that
%
\begin{eqnarray}\label{eq2.1}
0\leq E(U_{i1})=d_i ,\qquad\operatorname{var}(U_{i1}
)=1\qquad\mbox{for all
$i$} ,\qquad
\sup_{i\geq1} E(|U_{i1}|^3)<\infty ,
\end{eqnarray}
where $d_1,d_2,\ldots$ is a sequence of constants. The assumption that
$d_i\geq0$ is made here because, in the great majority of practical
applications, the hypothesis alternative to the null entails the zero
level being exceeded. Accordingly, the tests are one-sided, hence our
preoccupation with exceedences of a level. However, minor modifications
of our arguments permit the two-sided case to be treated.

Further, we assume that for an integer $\kappa >0$,
%
\begin{eqnarray}\label{eq2.2}
\begin{tabular}{@{}p{330pt}@{}}
the random vectors
$(U_{1i},U_{2i},\ldots)$, for $i\geq1$, are independent and
identically distributed,
the sequence of random variables $U_{11},U_{21},\ldots$ is
$\kappa $-dependent and
$\max_{i_1,i_2 \dvt i_1\neq i_2} \rho_{i_1i_2}<1$,
\end{tabular}
\end{eqnarray}
where $\rho_{i_1i_2}=\operatorname{corr}(U_{i_11},U_{i_21})$. The
third moment
condition in (\ref{eq2.1}) permits the variables $U_{ij}$ to have relatively
heavy-tailed distributions, for example, a Pareto distribution with
tail exponent greater than~3.

The assumption of short-range correlation in (\ref{eq2.2}) is, of course, an
oversimplification, but it reflects the low level of correlation that
is often observed in practice. For example, Messer and Arndt (\citeyear{2006Messer})
argue that correlation decays from about 0.08, at a separation of
approximately two base pairs, to about $0.01$ for a separation of ten
base pairs. Results reported by Mansilla \textit{et al.}~(\citeyear{2004Mansilla})
corroborate these figures if we assume that their data are normally
distributed. More generally, Almirantis and Provata (\citeyear{1999Almirantis}) give
evidence of both short-range and long-range correlation, depending on
the nature of the DNA or RNA under investigation.

The relationship between the group size, $n$, and the number of
hypothesis tests, $p$, is assumed to satisfy
%
\begin{equation}\label{eq2.3}
\log p=\mathrm{o}(n) .
\end{equation}
This allows the group size to be very much smaller than the number of tests.
In the absence of more detailed assumptions about the distributions of
the $U_{ij}$'s,
(\ref{eq2.3}) is necessary for the theorem we shall give below. To appreciate
why, note that if the $U_{ij}$'s are independent and identically
distributed with an atom at zero and, in particular, if $\delta \equiv
P(U_{ij}=0)>0$, then, with probability at least $\delta ^{n}$, the
$t$-statistic $T_1$ assumes the indeterminate value~$0/0$. In such
cases, we shall take $T_1=1$, but in order for the theorem to have a
meaningful interpretation when $t$ is the $(1-p^{-1})$-level quantile of
the standard normal distribution, it is essential that the probability
that $T_1=0/0$ be of smaller order than~$p^{-1}$. Therefore, we require
$\delta ^n=\mathrm{o}(p^{-1})$ for all $0<\delta <1$ and this assumption is
equivalent to~(\ref{eq2.3}).

Define
%
\begin{equation}\label{eq2.4}
\alpha =\frac{1}{4} \min_{i_1,i_2 \dvt i_1\neq i_2} (1-\rho
_{i_1i_2})
\end{equation}
and $\gamma =\alpha +1$. Condition (\ref{eq2.2}) implies that $\alpha >0$
and, of course,
$\alpha
\leq\frac{1}{2}$. Given $\eta>0$, let $t=t(p)$ satisfy
%
\begin{equation}\label{eq2.5}
(1+\eta) \sqrt{2 \gamma ^{-1} \log p} \leq t=\mathrm{O}\bigl(\sqrt{\log
p}\bigr) ,\qquad\max_{1\leq i\leq p} d_i=\mathrm{o}\bigl(t/\sqrt n\bigr) ,
\end{equation}
where $d_1,d_2,\ldots$ are as in~(\ref{eq2.1}). If $t$ satisfies the first part
of (\ref{eq2.5}), then any function $\phi$ which satisfies, as $t\to\infty $,
%
\begin{equation}\label{eq2.6}
\phi(t)=\exp\{\mathrm{o}(t^2)\}
\bigl\{\exp\bigl(-\tfrac{1}{4} t^2\bigr)+p \exp(-\gamma
t^2/2
)\bigr\}
\end{equation}
converges to zero as $p\to\infty $. In the arguments in Section~\ref{sec3}, we shall
use this notation generically; while $\phi$ will satisfy (\ref{eq2.6}), it will
alter from one appearance to another. Strictly speaking, it is not
essential to take $p$ to diverge. Although that condition motivates the
assumption of divergent $t$ and is, in turn, motivated by the
contemporary high-dimensional problems that led to this work, it is not
necessary for the theorem below.

\begin{theorem}\label{thm1} If $(\ref{eq2.2})$--$(\ref{eq2.6})$ hold, then there exists a
probability space on which are defined random variables
$T_1^{\operatorname{new}}
,\ldots
,T_p^{\operatorname{new}}$ and $T_1',\ldots,T_p'$ such that \textup{(i)}~the
joint distribution\vspace{1pt}
of $T_1^{\operatorname{new}},\ldots,T_p^{\operatorname{new}}$ is
identical to that of $T_1,\ldots,T_p$;
\textup{(ii)}~the random variables\vspace{1pt} $T_1',\ldots,T_p'$ are independent and
distributed, respectively, as $T_1,\ldots,T_p$; and \textup{(iii)}~with
probability equal to $1-\phi(t)$, the exceedences of $t$ by
$T_1^{\operatorname{new}}
,\ldots,T_p^{\operatorname{new}}$ occur at the same indices and take
the same\vspace{1pt} values as
the exceedences of $t$ by $T_1',\ldots,T_p'$.
\end{theorem}

To interpret the theorem, note that we would normally expect the
dependent data set $T_1,\ldots,T_p$ to exhibit clusters of level
exceedences, rather than the single, isolated exceedences associated
with the independent sequence $T_1',\ldots,T_p'$. The fact that the
$T_i$ process (or, equivalently,\vspace{1pt} the $T_i^{\operatorname{new}}$
process) behaves like
the $T_i'$ process in the case of large exceedences reflects the fact
that, since the marginal distribution of a $t$-statistic is relatively
light-tailed (if $n$ is sufficiently large -- see (\ref{eq2.3})), exceedences
of a high level are rare and so are unlikely to occur together. The
case of low-level exceedences is a very different matter, of course,
and so we would expect the theorem to fail if the lower bound for $t$,
in the first part of (\ref{eq2.5}), were relaxed too far.

\subsection{Applications}\label{sec2.2}

In this section, we treat the case of the null hypothesis, where
$d_i=0$ for each~$i$. This would be assumed in most applications of
Theorem~\ref{thm1} since it represents the setting that is conventionally used
for calibration.

The theorem implies that, in a strong sense, exceedences down to those
of the level $(1+\eta) (2 \gamma ^{-1} \log p)^{1/2}$ are
identical to the
ones that would occur in the case of independent tests. Now, the
probability associated with an exceedence of $(1+\eta) (2 \gamma
^{-1}
\log p)^{1/2}$ is, for small\vspace{1pt} $\eta$, approximately~$p^{-1/\gamma }$.
Therefore, false discoveries at probability levels of approximately
$p^{-1/\gamma }$, and at lower levels, can be adequately controlled by
assuming that the tests are\vspace{1pt} independent, even when they are not. Note
that $\gamma ^{-1}<1$ and that the false-discovery level controlled by the
conventional family-wise error rate is only~$p^{-1}$.

Next, we discuss the sorts of calculations that are enabled by
Theorem~\ref{thm1}. Let $Q_j$ denote the number of indices $i\in[1,p]$ for which
$T_i$ lies in the interval $(t_j,t_{j-1}]$, where $j\geq1$ and $t_j$ is
determined by $P(T_1>t_j)=j\beta /p$, with $\beta >0$ held fixed. (We take
$t_0=\infty$.) If $T_1,\ldots,T_p$ were fixed, then the joint
distribution of the $k+1$ random variables\vspace{1pt} $Q_1,\ldots,Q_k,p-\sum_j Q_j$
would be exactly multinomial with parameters $p$ and $q_1,\ldots
,q_k,1-\sum_j q_j$, where $q_j=P(T_1\in(t_j,t_{j-1}])$.
Theorem~\ref{thm1}
implies that for the dependent process $T_1,\ldots,T_p$, and for any
$k_0=k_0(p)$ for which $t_{k_0}$ satisfies (\ref{eq2.5}), any simultaneous
probability calculation based on the multinomial result, but applied to
the actual $T_j$ process rather than an idealized process with
independent marginals, is valid, provided that $k\leq k_0$ and the
final computed probability is quantified by adding an error which is
stated to be of order $\exp\{\mathrm{o}(t_k^2)\} \{\exp(-t_k^2/4)+p \exp
(-\gamma
t_k^2/2)\}$. The latter probability converges to zero, even if $k=k_0$
is taken as large as $p^{(\gamma _1-1)/\gamma _1}$, where $\gamma
_1\in(1,\gamma )$.

From this point, simultaneous multinomial probability calculations
based on $Q_1,\ldots,Q_k$, familiar from the well-understood case of
independent test statistics, can be used to construct rules for
controlling FWER or false discovery rate~(FDR); see, for example,
Benjamini and Hochberg~(\citeyear{1995Benjamini}). Wang and Hall (\citeyear{2009Wang}) have shown that,
under the assumption of finite third moments, highly accurate
approximations are available for the marginal distribution of~$T_1$.
Such calculations, which justify standard normal, Student's $t$- or
bootstrap approximations to the marginal distribution of $T_1$, are
already widely used in practice (see Section~\ref{sec1}), in conjunction with
the independence assumption, when controlling false discovery rates.
Our paper provides justification for these methods.

More generally, Theorem~\ref{thm1} implies that if a probability statement about
what the process $T_1,\ldots,T_p$ does above the level $t$ is founded
on the assumption of independence, then, no matter how complex or
convoluted the statement might be, the claimed probability level is
accurate to within~$\phi(t)$.

To give an example of calculations based on Theorem~\ref{thm1}, take $p\leq
p_0=10^6$, $n=100$ and $t=5.052$, the latter denoting the upper
$(1-p_0^{-1})$-level quantile of Student's $t$-distribution with $n-1=99$
degrees of freedom. Reflecting empirical evidence given in Section~\ref{sec2.1},
take $\gamma =\frac{1}{4} (1-0.1)+1=1.225$. Then, (\ref{eq2.5}) is in order; the
probability that at least one value of $p$ independent $t$ statistics,
each on 99 degrees of freedom, exceeds $t=5.052$ equals $0.010$,
$0.095$ and $0.63$ for $p=10^4$, $10^5$ and $10^6$, respectively; and
(\ref{eq2.6}) suggests that the errors in these levels are in error by less
than $30\%$, $20\%$ and $0.25\%$, respectively. Most likely, the errors
are much less than these since the asymptotic bound is derived only as
an upper bound. If we were to make a general probability statement
about exceedences of the level $5.052$ by the stochastic process of $t$
statistics, under the assumption of independence, then, despite the
process actually being $\kappa $-dependent rather than independent, we
would expect to make errors no greater than these respective values. In
the same general setting, relative error decreases to zero as $p$ and
$t$ increase. For example, in cases where $t$ solves $1-\Phi(t)\asymp
p^{-1}$, with $\Phi$ denoting the standard normal distribution function,
we have $\{\exp(-\frac{1}{4} t^2)+p \exp(-\frac{1}{2} \gamma
t^2)\}
/(1-p)^p=\mathrm{O}[\exp
(-\frac{1}{4} t^2)+t \exp\{-\frac{1}{2} (\gamma -1) t^2\}]\to
0$ as $t\to\infty $.

\subsection{Generalizations}\label{sec2.3}

Theorem~\ref{thm1} can be extended to other settings, in particular, to those
where (a)~a wider range of dependence, obtained by allowing $\kappa $ in
(\ref{eq2.2}) to diverge with $p$, is allowed; (b)~the value of $n$ for the
$i$th group equals $n_i$, depending on $i$, and (\ref{eq2.3}) is altered by
requiring\vspace{1pt} that $\log p=\mathrm{o}(\min_{i\leq p} n_i)$; (c)~weights $w_{ij}$ are
incorporated into the construction of the $t$-statistics $T_i$, by
defining $\bar{U}_i= n_i^{-1} \sum_j w_{ij} U_{ij}$, $S_i^2=
n_i^{-1} \sum_j
w_{ij}
^2 U_{ij}^2-\bar{U}_i^2$ and, as before, $T_i= n_i^{1/2} \bar{U}_i/S_i$.
Provided the weights satisfy
\[
\sup_{i,j} |w_{ij}|\leq C_1 ,\qquad
\inf_{1\leq i\leq p} n_i^{-1} \#\{j\dvt |w_{ij}|\geq C_2\}\geq C_3 ,
\]
where $C_1,C_2,C_3$ are positive constants not depending on $p$, the
proof in this more general case is as in Section~\ref{sec3}. However, the
statement of the theorem is then less elegant and less transparent, so
we do not give the more general version here. Incorporation of the
weights $w_{ij}$ permits the scope of the example above to be extended to
hypothesis-testing problems involving linear regression.\looseness=1

To indicate the types of results that can be achieved under longer
ranges of dependence, we shall discuss the case of a moving average,
\[
U_{ij}=\kappa ^{-1/2} \sum_{k=1}^\kappa  \varepsilon _{j,i+k} ,
\]
where $\kappa =\kappa (p)$ is permitted to diverge to infinity at a
rate not
exceeding $\log p$ and the independent disturbances $\varepsilon
_{ji}$ are all
distributed as $\varepsilon $, for which $E(\varepsilon )=0$ and
$E|\varepsilon |^3<\infty$. In
this setting, (\ref{eq2.2}) holds. We strengthen (\ref{eq2.3}) by asking that $\log
p=\mathrm{O}(n^{1/3})$. The definition of $t$ implicit in~(\ref{eq2.5}) can now be
refined to
\[
t=\sqrt{2 \gamma ^{-1} (\log p+A \log\log p)} ,
\]
where $A>0$ denotes a sufficiently large absolute constant. The
conclusions of Theorem~\ref{thm1} continue to hold, with a similar proof if we
replace $1-\phi(t)$ by $1-\mathrm{o}(1)$.

\section{\texorpdfstring{Proof of Theorem \protect\ref{thm1}}{Proof of Theorem 1}}\label{sec3}

\subsection{Step~1: Preliminaries}\label{sec3.1}

 The notation $D_1,D_2,\ldots$
will denote constants not depending on $n$ or~$p$. Let $Q_i=
n^{-1} \sum_j (U_{ij}-d_i)^2$, $R_i= n^{1/2} \bar{U}_i/Q_i^{1/2}$ and
note that
%
\begin{eqnarray}\label{eq3.1}
\begin{tabular}{@{}p{330pt}@{}}
for each $t>0$, the events $T_i>t$ and $R_i>t/(1+n^{-1}
t^2)^{1/2}$ are identical.
\end{tabular}
\end{eqnarray}
Also, note that
$R_i=(\sum_j V_{ij}+nd_i)/(\sum_j V_{ij}^2)^{1/2}$, where
$V_{ij}=U_{ij}-d_i$, $1\le j\le n$,\vspace{1pt} are independent and identically
distributed random variables
satisfying
\renewcommand{\theequation}{\arabic{section}.\arabic{equation}a}
\setcounter{equation}{1}
\begin{eqnarray}\label{eq3.2a}
0\leq E(V_{i1})=0 ,\qquad E(V_{i1}^2)=1\qquad\mbox{for
all $i$} ,\qquad
\sup_{i\geq1} E(|V_{i1}|^3)<\infty ;
\end{eqnarray}
cf.~(\ref{eq2.1}).

\subsection{Step~2: Probabilities of exceedences in ones and
twos}\label{sec3.2}

Using results of Wang and Hall (\citeyear{2009Wang}) (see also Wang (\citeyear{2005Wang})), it can be
shown that, for constants $D_1,D_2,D_3>0$, and whenever $0<s<D_1
n^{1/2}$,
\renewcommand{\theequation}{\arabic{section}.\arabic{equation}}
\setcounter{equation}{1}
\begin{equation}\label{eq3.2}
P(R_i>s)\leq D_2 s^{-1} \exp\bigl(D_3 s^3 n^{-1/2}-\tfrac{1}{2}
s^2+\sqrt
n d_i  s\bigr).
\end{equation}

\noindent
We also wish to prove the following related result for pairs of exceedences.

\begin{lemma*} Assume the conditions of Theorem~\ref{thm1}. There then exist
$D_4,D_5>0$ such that for all $i_1,i_2$ with $i_1\neq i_2$, and for all
$0<s<D_4 n^{1/2}$,
we have
\begin{eqnarray}\label{eq3.3}
&&\sup_{1\leq|i_1-i_2|\leq k_2-k_1}
P(R_{i_1}>s,R_{i_2}>s)\nonumber
\\[-8pt]\\[-8pt]
&&\quad \leq5 \exp\bigl\{-\tfrac{1}{2}
(1+\alpha ) s^2
+D_5 n^{-1/2} s^3 +2  \sqrt
n (d_{i_1}+d_{i_2}) s\bigr\},\nonumber
\end{eqnarray}
where $\alpha$ is as in~$(\ref{eq2.4})$.
\end{lemma*}

To establish the lemma, we write
\begin{eqnarray*}
 U_{i_11}^{(1)}&=&V_{i_11} I(|V_{i_11}|\leq
n^{1/2}/s ,  |V_{i_21}|\leq n^{1/2}/s) ,\qquad
U_{i_11}^{(2)}=V_{i_11}-U_{i_11}^{(1)} ,
\\[2pt]
U_{i_21}^{(1)}&=&V_{i_21} I(|V_{i_11}|\leq n^{1/2}/s ,
|V_{i_21}|\leq n^{1/2}/s) , \qquad
U_{i_21}^{(2)}=V_{i_21}-U_{i_11}^{(1)} .
\end{eqnarray*}
By virtue of (\ref{eq3.2a}), simple calculations show that
\begin{eqnarray*}
\bigl|E\bigl(U_{i_11}^{(1)}+U_{i_21}^{(1)}\bigr)\bigr|
&\leq&D_6 s^2/n ,
\\[3pt]
E\bigl\{\bigl(U_{i_11}^{(1)}\bigr)^2+\bigl(U_{i_21}^{(1)}\bigr)^2\bigr\}
&=&2+\mathrm{O}(1)  s/\sqrt n ,
\\[3pt]
E\bigl(U_{i_11}^{(1)}+U_{i_21}^{(1)}\bigr)^2
&=&E\bigl\{(V_{i_11}+ V_{i_21})-\bigl(U_{i_11}^{(2)}+U_{i_21}^{(2)}\bigr)\bigr\}^2
\\[3pt]
&\leq& 2 (1+\rho_{ij})+D_6 s/\sqrt n ,
\\[3pt]
E\bigl|U_{i_11}^{(1)}+U_{i_21}^{(1)}\bigr|^3&\leq& D_6 .
\end{eqnarray*}
These results, and the bound $\mathrm{e}^x\leq1+x+\frac{1}{2} x^2+\frac
{1}{6} |x|^3 \mathrm{e}^x$
(valid for all real $x$), imply that, with $h=s/\sqrt{n} $,
%
\begin{eqnarray}\label{eq3.4}
&&E\biggl[\exp\biggl\{\frac{1}{2} h (V_{i_11}+V_{i_21})
-\frac{1}{4} h^2 (V_{i_11}^2+V_{i_21}^2)\biggr\}
I(|V_{i_11}|\leq n^{1/2}/s , |V_{i_21}|\leq
n^{1/2}/s)\biggr]\nonumber
\\[-7pt]\\[-7pt]
&&\quad
\leq1+(\rho_{ij}-1) {s^2\over4n}+ D_7 {s^3\over
n^{3/2}} ,\nonumber
\\[2pt] \label{eq3.5}
&&E\biggl[\exp\biggl\{\frac{1}{2} h (V_{i_11}+V_{i_21})
-\frac{1}{4} h^2 (V_{i_11}^2+V_{i_21}^2)\biggr\}
I(|V_{i_11}|\geq n^{1/2}/s ,\mbox{ or }|V_{i_21}|\geq
n^{1/2}/s)\biggr]\nonumber
\\[2pt]
&&\quad\leq
\exp \{P(|V_{i_11}|>n^{1/2}/s)
+P(|V_{i_21}|>n^{1/2}/s)\}
\\[2pt]
&&\quad\leq D_8 \bigl(s/\sqrt n\bigr)^3 .\nonumber
\end{eqnarray}

\noindent
Results (\ref{eq3.4}) and (\ref{eq3.5}), together with the independence of
$V_{i_1k}$ for each $i_1$, imply that, for $s\leq D_9\sqrt n$,
with $D_9$ sufficiently small,
\begin{eqnarray}\label{eq3.6}
 &&E\biggl[\exp\biggl\{{h \over2} \sum_k(
V_{i_1k}+V_{i_2k})- {h^2\over4} \sum_k(
V_{i_1k}^2+V_{i_2k}^2)\biggr\}\biggr]\nonumber
\\[-8pt]\\[-8pt]
&&\quad \leq\biggl\{1+(\rho_{ij}-1) {s^2\over4n} +D_{10} {s^2\over
n^{3/2}}\biggr\}^{\!n}
\leq\exp\biggl\{(\rho_{ij}-1) {s^2\over4}+D_{11}  {s^3\over\sqrt
n}\biggr\} .\nonumber
\end{eqnarray}

Define $\varepsilon ^2=(1-\rho_{ij})/8$. It follows from (\ref{eq3.6}) that
whenever $s\leq D_9 \sqrt n$, with $D_9$ sufficiently small,
we have
%
\begin{eqnarray}\label{eq3.7}
 \pi_{1n} &\equiv& P\biggl\{2
h \sum_k(V_{i_1k}+V_{i_2k})-h^2
\sum_k(V_{i_1k}^2+V_{i_2k}^2)\nonumber
\\[2pt]
&&{}\quad\  +2 \sqrt n (d_{i_1}+d_{i_2}) s
\geq2 s^2 (1-\varepsilon ^2)\biggr\}\nonumber
\\[2pt]
&\leq&\exp\biggl\{2 \sqrt n (d_{i_1}+d_{i_2}) s-\frac{1}{2}s^2
(1-\varepsilon^2)\biggr\}
\\[2pt]
&&{}\times
E\biggl[\exp\biggl\{\frac{1}{2} h \sum_k(V_{i_1k}+V_{i_2k})
-\frac{1}{4} h^4 \sum_k(V_{i_1k}^2+V_{i_2k}^2)\biggr\}\nonumber
\biggr]
\\[2pt]
&\leq& \exp\biggl\{-\frac{1}{2} (1+\alpha ) s^2 +D_{12} {s^3\over
\sqrt
n}+2  \sqrt n (d_{i_1}+d_{i_2}) s\biggr\} ,\nonumber
\end{eqnarray}
where $\alpha $ is as defined in~(\ref{eq2.4}). Write $\Omega
_n=(1-\varepsilon ,1+\varepsilon )$ and
note that if $0<\varepsilon <\frac{1}{2}$, then
\begin{eqnarray*}
&& \biggl\{\sum_kV_{i_2k}+n  d_{i_2} \geq
s (n Q_{i_2}^2)^{1/2} , Q_{i_2}\in\Omega _n \biggr\}
\\[2pt]
&&\quad \subseteq\biggl\{2 h \sum_kV_{i_2k}-h^2 \sum_k
V_{i_2k}^2+2 \sqrt n d_{i_2} s \geq
s^2 (1-\varepsilon ^2)\biggr\},
\end{eqnarray*}
where $h=s/\sqrt n$. It can be shown that
%
\begin{eqnarray}\label{eq3.8}
 &&P(R_{i_1}>s,R_{i_2}>s)\nonumber
 \\[2pt]
 &&\quad =P\biggl\{\sum_k
V_{i_1k}+n  d_{i_1} \geq
s (n Q_{i_1}^2)^{1/2} ,
\sum_k
V_{i_2k}+n  d_{i_2}\geq
s (n Q_{i_2}^2)^{1/2}\biggr\}
\\[2pt]
&&\quad \leq\pi_{1n}+\pi_{2n}+\pi_{3n}+\pi_{4n}+\pi_{5n}
,\nonumber
\end{eqnarray}

\noindent
where
\begin{eqnarray*}
 \pi_{2n}&=&P\biggl\{\sum_kV_{i_1k}+n  d_{i_1} \geq
s (n Q_{i_1}^2)^{1/2} ,  Q_{i_1}\geq1+\varepsilon
\biggr\},
\\
\pi_{3n}&=&P\biggl\{\sum_kV_{i_2k}+n  d_{i_2} \geq
s (n Q_{i_2}^2)^{1/2} ,  Q_{i_2}\geq1+\varepsilon
\biggr\},
\\
\pi_{4n}&=&P\biggl\{\sum_kV_{i_1k}+n  d_{i_1} \geq
s (n Q_{i_1}^2)^{1/2} ,  Q_{i_1}\leq1-\varepsilon
\biggr\} ,
\\
\pi_{5n} &=&P\biggl\{\sum_kV_{i_2k}+n  d_{i_2} \geq
s (n Q_{i_2}^2)^{1/2} ,  Q_{i_2}\leq1-\varepsilon
\biggr\}.
\end{eqnarray*}
Property (\ref{eq3.3}) will follow from (\ref{eq3.7}) and (\ref{eq3.8}) if we prove that
there exists $D_{13}>0$ such that, for $s\leq D_{13} n^{1/2}$,
%
\begin{eqnarray}\label{eq3.9}
\pi_{kn}\leq\exp\biggl\{-\frac{1}{2} (1+\alpha ) s^2+D_{12}
{s^3\over\sqrt
n}  +2 \sqrt n (d_{i_1}+d_{i_2}) s\biggr\}
\end{eqnarray}
for $k=2,3,4,5$.

Our proof of (\ref{eq3.9}) is based on arguments of Shao (\citeyear{1999Shao}) (see also
the proof of Proposition~4.2 of Wang and Hall (\citeyear{2009Wang})) and uses the
following result: if $EX=0$, $EX^2=1$ and $E|X|^3<\infty$, then for
any $\lambda >0$, $\theta >0$ and $x>0$,
%
\begin{eqnarray}\label{eq3.10}
E[\exp\{\lambda  b X-\theta  (bX)^2\}]
=1+(\lambda ^2-\theta ) n^{-1} x^2
+A(\lambda ,\theta ) n^{-3/2} x^3 E|X|^3 ,
\end{eqnarray}
where $b=x/\sqrt n$ and $A(\lambda , \theta )$ depends only on
$\lambda $ and
$\theta $. This result is a special case of Lemma~1 of Shao (\citeyear{1999Shao}). Also,
note that
\begin{eqnarray*}
\pi_{2n}&\leq& \pi_{2n}^{(1)} +P\biggl\{\sum_kV_{i_1k} +n d_{i_1}\geq
s (n Q_{i_1}^2)^{1/2} ,  Q_{i_1}\ge3\biggr\}
\\
&\leq& \pi_{2n}^{(1)}+\pi_{2n}^{(2)}+\pi_{2n}^{(3)} ,
\end{eqnarray*}
where, noting that $\sqrt n |d_i|\leq s/5$, we define
\begin{eqnarray*}
 \pi_{2n}^{(1)}&=&P\biggl\{\sum_kV_{i_1k} +n d_{i_1}\geq
s (n Q_{i_1}^2)^{1/2} ,   1+\varepsilon \leq
Q_{i_1}<3\biggr\}
,\\
\pi_{2n}^{(2)}&=&P\biggl\{\sum_k
V_{i_1k} I(|V_{i_1k}|> n^{1/2}/s) \geq s \biggl(\sum_k
V_{i_1k}^2\biggr)^{\!1/2} \biggr\} ,
\\
\pi_{2n}^{(3)}&=&P\biggl\{
\sum_k
V_{i_1k} I(|V_{i_1k}|\leq n^{1/2}/s) \geq3 s \sqrt
n/2\biggr\} .
\end{eqnarray*}
If the random variable $H$ has the $\operatorname{Bi}(n,p)$ distribution and if
$a>0$, then $P(H>an)\leq(\mathrm{e}p/a)^{an}$ and so
\begin{eqnarray*}
 \pi_{2n}^{(2)}&\leq& P\biggl\{\sum_k
I(|V_{i_1k}|> n^{1/2}/s) \geq s^2\biggr\}
\\
&\leq&\{s^{-2} 12 n P(|V_{i_1k}|>n^{1/2}/s)
\}
^{\!s^2}
\leq\frac{1}{2} \mathrm{e}^{-s^2}
\end{eqnarray*}
for $s\le D_{14} \sqrt n$, with $D_{14}$ sufficiently small.
Arguments similar to those in the proof of (\ref{eq3.7}) yield that
$\pi_{2n}^{(3)}\leq\frac{1}{2} \mathrm{e}^{-s^2}$ for $s\le D_{14} \sqrt n$ with
$D_{14}$ sufficiently small. To estimate $\pi_{2n}^{(1)}$, we write
$\mathcal{S}_1=\{(x,y)\dvt x\ge s \sqrt y, s^2 (1+\varepsilon )^2\le
y\le9 s^9\}$.
It follows from (\ref{eq3.10}) with $\lambda =1$, $\theta =\frac{1}{6}$ and
$X=V_{i_11}$ that,
with $h=s/\sqrt n $,
\begin{eqnarray*}
 \pi_{2n}^{(1)}&=&P\biggl\{\biggl(h \sum_kV_{i_1k}+\sqrt
n
d_i s,
 h^2 \sum_kV_{i_1k}^2\biggr)\in\mathcal{S}_1\biggr\}
 \\
 &\leq&
E\biggl[\exp\biggl(h \sum_kV_{i_1k}-\frac{1}{6} h^2 \sum_k
V_{i_1k}^2+\sqrt n d_i s\biggr)
 \exp\Bigl\{-\inf_{(x, y)\in
\mathcal{S}_1} (x-y/6)\Bigr\}\biggr]
\\
&\leq& \exp\biggl\{\biggl(\frac{1}{2}-\frac{1}{6}\biggr) s^2
-s^2 (1+\varepsilon )+\frac{1}{6} s^2 (1+\varepsilon )^2+\sqrt n
d_i s+D_{15} s^3
n^{-1/2}
\biggr\}
\\
&\leq& \exp\biggl\{-\frac{1}{2}  s^2-(5 \varepsilon  s^2
/8)
+\sqrt n d_i s+\bigl(D_{15}  s^3/\sqrt n\bigr)\biggr\}
\\
&\leq& \exp\biggl\{-\frac{1}{2} (1+\alpha ) s^2+\sqrt n d_i s
+\bigl(D_{15}
s^3/\sqrt
n\bigr)\biggr\} ,
\end{eqnarray*}
where we have used the fact that the function $f(y)=s \sqrt
y-\frac{1}{6} y$ is increasing in $s^2 (1+\varepsilon )^2\leq y\leq
9 s^2$.
Combining all of the above estimates, we obtain
\[
\pi_{2n}\leq\exp\bigl\{-\tfrac{1}{2} (1+\alpha ) s^2+\sqrt n
d_i s
+\bigl(D_{15} s^3 /\sqrt n\bigr)\bigr\} .
\]
Similarly, we may prove (\ref{eq3.9}) for $k=3$.

Put $\mathcal{S}_2=\{(x,y)\dvt x\ge s \sqrt y,  y\le(1-\varepsilon
)^2 s^2\}$. It
follows from (\ref{eq3.10}) with $\lambda =1$, $\theta =2$ and $X=V_{i_11}$
that, with
$h=s/\sqrt n$,
\begin{eqnarray*}
 \pi_{4n}&=&P\biggl\{\biggl(h \sum_kV_{i_1k}+\sqrt n
d_i
s
,
h^2 \sum_kV_{i_1k}^2\biggr)\in\mathcal{S}_2\biggr\}
\\
&\leq& E\biggl[\exp\biggl(h \sum_kV_{i_1k} -2 h^2 \sum_k
V_{i_1k}^2+\sqrt n d_i s\biggr)
 \exp\Bigl\{-\inf_{(x,y)\in
\mathcal{S}_2} (x-2y)\Bigr\}\biggr]
\\
&\leq& \exp\bigl\{-1.5
s^2-s^2
(1-\varepsilon )
+2 s^2(1-\varepsilon )^2+\sqrt n d_i s+\bigl(D_{16} s^3/\sqrt
n\bigr)\bigr\}
\\
&\leq& \exp\biggl\{-\frac{1}{2} s^2-2 \varepsilon  s^2+\sqrt n
d_i s
+\bigl(D_{16} s^3/\sqrt n\bigr)\biggr\}
\\
&\leq& \exp\biggl\{-\frac{1}{2} (1+\alpha ) s^2+\sqrt n d_i s
+\bigl(D_{16}
s^3/\sqrt n\bigr)\biggr\} .
\end{eqnarray*}
Similarly, we may prove (\ref{eq3.9}) for $k=5$. This completes the
derivation of (\ref{eq3.9}) and, hence, also the proof of the lemma.

\subsection{Step~3: Blocks and expected numbers of level
exceedences}\label{sec3.3}

Partition the set of positive integers into small blocks, each of
length $\kappa +1$, where $\kappa $ is as in (\ref{eq2.2}), and large blocks,
each of
length $\ell$, where $\ell$ is a divergent function of~$p$. We shall take
%
\begin{eqnarray}\label{eq3.11}
\ell\sim\exp\bigl(\tfrac{1}{4} s^2\bigr) ,
\end{eqnarray}
where $s\to\infty $ as $p$ increases. The integers in each block are
consecutive, each consecutive pair of large blocks is separated by a
small block and the block furthest to the left is a large block. Let
the small blocks be $b_1,b_2,\ldots$ and the large blocks be
$B_1,B_2,\ldots,$ indexed such that the order of the blocks is
$B_1,b_1,B_2,b_2,\ldots.$ Let $B=B_1=\{1,\ldots,\ell\}$ denote the
first large block and let $N_1$ be the number of indices $i\in B$ for
which $R_i>s$. We wish to develop a bound for $E\{N_1 I(N_1\geq2)\}$.
Identical bounds can be derived, uniformly in the block indices, for
the versions of $N_1$ in the case of blocks $B_2,B_3,\ldots;$ for
notational simplicity, we focus solely on $B_1$.

By H\"{o}lder's inequality,
%
\begin{eqnarray}\label{eq3.12}
E\{N_1 I(N_1\geq2)\}
\leq(EN_1^{a_1})^{1/a_1} P(N_1\geq2)^{1/a_2} ,
\end{eqnarray}
where $a_1,a_2>1$ satisfy $a_1^{-1}+a_2^{-1}=1$. Define $d^0=\sqrt n
\max
_{1\leq i\leq p} d_i$. In view of (\ref{eq3.2}) and (\ref{eq3.3}),
%
\begin{eqnarray}\label{eq3.13}
P(N_1\geq2)
&=&P(\mbox{for some $i_1,i_2\in B$ with $i_1<i_2$,
$R_{i_1},R_{i_2}>s$})\nonumber
\\
&\leq&\sum_{i_1=1}^{\ell-1} \sum_{i_2=i_1+1}^\ell
P(R_{i_1}>s,R_{i_2}>s)\nonumber
\\
&=&\sum_{i_1=1}^{\ell-1} \sum_{i_2=i_1+1}^{\min(i_1+\kappa +1,\ell
)} P(R_{i_1}>s,R_{i_2}>s)
\\
&&{}+\sum_{i_1=1}^{\ell-1} \sum_{i_2=\min(i_1+\kappa +2,\ell)}^\ell
P(R_{i_1}>s) P(R_{i_2}>s)\nonumber
\\
&\leq& D_{17} \exp(D_{18} s^3 n^{-1/2}+D_{19} d^0 s)\nonumber
\\
&&{}\times
\biggl[\ell \exp\biggl\{-\frac{1}{2} (\alpha +1) s^2\biggr\}
+\ell^2 \exp(-s^2)\biggr] .\nonumber
\end{eqnarray}
Noting that $N_1$ can be written as $\kappa +1$ sums of $\ell
/(\kappa +1)$
independent and identically distributed random variables and using
calculations based on the binomial distribution, it can be shown that,
for the choice of $\ell$ at (\ref{eq3.11}), $E(N_1^{a_1})$ is bounded as
$p\to\infty
$ for each $a_1>0$. Hence, using (\ref{eq3.12}) and (\ref{eq3.13}), we deduce that for
each $\eta_2\in(0,1)$,
\begin{eqnarray}\label{eq3.14}
E\{N_1 I(N_1\geq2)\}&\leq& D_{20} \exp(D_{21} s^3 n^{-1/2}
+D_{22}
d^0 s)\nonumber
\\[-8pt]\\[-8pt]
&&{}\times
\bigl[\ell \exp\bigl\{-\tfrac{1}{2} (\alpha +1) s^2\bigr\}+\ell
^2 \exp
(-s^2
)\bigr]^{1-\eta_2} .\nonumber
\end{eqnarray}

Write $N_2$ for the number of exceedences of $s$ that occur in the
union of the small blocks $b_j$ that intersect the interval $[1,p]$.
There are $\mathrm{O}(p/\ell)$ such small blocks and each is of length $\kappa +1$,
so, by (\ref{eq3.2}),
%
\begin{eqnarray}\label{eq3.15}
E(N_2)\leq D_{23} p \ell^{-1} P(R_1>s)
\leq D_{24} p \ell^{-1} \exp\bigl(D_3 s^3 n^{-1/2}-\tfrac{1}{2}
s^2+d^0
s\bigr) .
\end{eqnarray}
Provided we choose $s=s(p)$ to diverge to infinity in such a manner that
%
\begin{eqnarray}\label{eq3.16}
s=\mathrm{O}\bigl(\sqrt{\log p}\bigr) ,\qquad d^0=\mathrm{o}(s) ,
\end{eqnarray}
it follows from (\ref{eq2.3}) that $s^3 n^{-1/2}+d^0 s=\mathrm{o}(s^2)$ and so (\ref{eq3.14})
entails that
\begin{eqnarray*}
E\{N_1 I(N_1\geq2)\}
=\exp\{\mathrm{o}(s^2)\}
\bigl[\ell \exp\bigl\{-\tfrac{1}{2} (\alpha +1) s^2\bigr\}
+\ell^2 \exp(-s^2)\bigr]^{1-\eta_2} .
\end{eqnarray*}
Since this is true for each $\eta_2>0$, we have
%
\begin{eqnarray}\label{eq3.17}
E\{N_1 I(N_1\geq2)\}
=\exp\{\mathrm{o}(s^2)\}
\bigl[\ell \exp\bigl\{-\tfrac{1}{2} (\alpha +1) s^2\bigr\}
+\ell^2 \exp(-s^2)\bigr] .
\end{eqnarray}

\subsection{Step~4: Bound for $P(N_1\geq1)$, and related
bounds}\label{sec3.4}

Let $N_3$ denote the number of exceedences of $s$ which come from large
blocks $B_j$, $1\leq j\leq m$, that have two or more exceedences. Write
$\sum_j \pi_j$ for the sum over $1\leq j\leq m$ of the probability
$\pi
_j$ that $R_i>s$ for some~$i\in B_j$. Then (a)~the expected number of
exceedences of $s$ by $R_1,\ldots,R_p$ equals $\sum_{i\leq p}
P(R_i>s)$ and is less than or equal to $\sum_j \pi_j+E(N_2)+E(N_3)$;
(b)~the\vspace{1pt} expected number of exceedences in (a) is greater than or equal
to $\sum_{j\leq m-1} \pi_j$; and (c)~since $P(N_1\geq1)\leq
E(N_1)=\ell
 P(R_1>s)$ and $P(R_1>s)$ satisfies (\ref{eq3.2}), we have
%
\begin{eqnarray}\label{eq3.18}
\pi_1=P(N_1\geq1)\leq P^0\equiv D_2 s^{-1} \ell
\exp\bigl(D_3 s^3 n^{-1/2}-\tfrac{1}{2} s^2+d^0 s\bigr)
\end{eqnarray}
and an identical bound holds for $\pi_1,\ldots,\pi_m$, in particular,
(d)~$\pi_m\leq P^0$. Results (a)--(d) imply that
%
\begin{eqnarray}\label{eq3.19}
\Biggl|\sum_{j=1}^m\pi_j-\sum_{i=1}^pP(R_i>s)\Biggr|
\leq E(N_2)+E(N_3)+P^0 .
\end{eqnarray}
Since $E(N_3)\leq m E\{N_1 I(N_1\geq2)\}$, $m=\mathrm{O}(p/\ell)$ and bounds
for $E\{N_1 I(N_1\geq2)\}$, $E(N_2)$ and $P(N_1\geq1)$ are given by
(\ref{eq3.17}), (\ref{eq3.15}) and (\ref{eq3.18}), it follows that (\ref{eq3.19}) entails, on taking
$\ell$ as in~(\ref{eq3.11}),
\begin{eqnarray}\label{eq3.20}
\Biggl|\sum_{j=1}^m\pi_j-\sum_{i=1}^pP(R_i>s)\Biggr|
 =\exp\biggl\{-\frac{1}{4} s^2+\mathrm{o}(s^2)\biggr\}
\biggl\{1+p \exp\biggl(-\frac{1}{4} s^2-\frac{1}{2}
\alpha  s^2\biggr)\biggr\} .\quad
\end{eqnarray}

\subsection{Step~5: Probabilities of level exceedences}\label{sec3.5}

Let $\mathcal{F}$ denote the event that (a)~there are no exceedences
of $s$ in
any of the small blocks that are wholly contained within $[1,p]$;
(b)~in each of the large blocks that is wholly contained within
$[1,p]$, there is at most one exceedence of~$s$; and (c)~there are no
exceedences of $s$ in any block fragment that overlaps the end
point~$p$. Write $\mathcal{G}$ for the complement of~$\mathcal{F}$.
Results (\ref{eq3.15}),
(\ref{eq3.17}) and (\ref{eq3.18}) imply that, with $\ell$ given by (\ref{eq3.11}) and assuming
that (\ref{eq3.16}) holds,
%
\begin{eqnarray}\label{eq3.21}
P(\mathcal{G})
\leq\exp\bigl\{-\tfrac{1}{4} s^2+\mathrm{o}(s^2)\bigr\}
\bigl\{1+p \exp\bigl(-\tfrac{1}{4} s^2-\tfrac{1}{2} \alpha
s^2\bigr)\bigr\} .
\end{eqnarray}
Therefore, in order for $P(\mathcal{G})\to 0$, it is sufficient that
for some
$\eta_3\in(0,1)$ and all sufficiently large $p$,
we have
%
\begin{eqnarray}\label{eq3.22}
(1+\eta_3) \sqrt{2 \gamma ^{-1} \log p}
\leq s=\mathrm{O}\bigl(\sqrt{\log p}\bigr) ,
\end{eqnarray}
where $\gamma $ is as defined in Section~\ref{sec2}. This choice of $s$
satisfies~(\ref{eq3.16}) and so if $s$ is given by (\ref{eq3.22}), then $P(\mathcal{G})$
satisfies~(\ref{eq3.21}).

\subsection{Step~6: Strong approximation}\label{sec3.6}

Let $M_j$, $1\leq j\leq m$, be the number of times that $R_i>s$ for
$i\in B_j$. Then, the number, $N$, say, of blocks $B_j$ for which
$M_j\geq1$ is distributed as $\sum_j I_j$, where the random variables
$I_j$ are independent, $I_j=1$ if $M_j\geq1$ and $I_j=0$ otherwise. As
before, we define $\pi_j=P(M_j\geq1)$. Conditional on $N$ and on the
events ``$M_{j_1}\geq1$'' and ``$M_{j_2}\geq1$,'' where $1\leq
j_1<j_2\leq m$, the sequences $\{R_i\dvt i\in B_{j_1}\}$ and $\{R_i\dvt i\in
B_{j_2}\}$ are independent.

Order the blocks $B_j$ for which $M_j\geq1$, giving $B_{J_1},\ldots
,B_{J_N}$, where $1\leq J_1<\cdots<\break J_N\leq m$, and let $W_k$ denote a
value of $R_i$ for which $R_i>s$, randomly chosen among such values for
which $i\in B_{J_k}$. Write $i=I_k$ for the index of the value of $R_i$
that is chosen as~$W_k$. Then, conditional on $N$, the random variables
$W_1,\ldots,W_N$ are independent and identically distributed as~$R(s)$,
$J_1,\ldots,J_N$ is a set of integers chosen independently and randomly
from $1,\ldots,m$ and $I_k$ is uniformly distributed among indices
in~$B_{J_k}$.

Let $R_1',\ldots,R_p'$ be independent random variables having the
distributions of~$R_1,\ldots,\break R_p$, respectively, let $M_j'$ denote
the number of times that $R_i'$ exceeds $s$ for $i\in B_j$ and put $\pi
_j'=P(M_j'\geq1)$. The numbers $N'$ of blocks $B_j$ for which
$M_j'\geq
1$ are distributed as $\sum_j I_j'$, where the random variables $I_j'$
are independent and $I_j'=1$ if $M_j'\geq1$, $I_j'=0$ otherwise. An
argument similar to, but simpler than, that leading to (\ref{eq3.20}) shows that
%
\begin{eqnarray}\label{eq3.23}
\sum_{j=1}^m|\pi_j-\pi_j'|
\leq\exp\biggl\{-\frac{1}{4} s^2+\mathrm{o}(s^2)\biggr\}
\biggl\{1+p \exp\biggl(-\frac{1}{4} s^2-\frac{1}{2} \alpha
s^2\biggr)\biggr\} .
\end{eqnarray}

By enlarging the probability space if necessary, we can think of $N$ as
denoting the number out of $m$ independent and random variables
$U_1,\ldots,U_j$, each uniformly distributed on $[0,1]$, which lie in
the respective intervals~$[0,\pi_j]$. Take $N'$ to be the number of
$U_i$'s that lie in~$[0,\pi_j']$. Then,
%
\begin{eqnarray}\label{eq3.24}
P(N=N')\geq1-\sum_{j=1}^m|\pi_j-\pi_j'| .
\end{eqnarray}

We have already constructed sequences $W_1,\ldots,W_N$, $I_1,\ldots
,I_N$ and $J_1,\ldots,J_N$. If $N'>N$, then, conditional on these
quantities and on $N$ and $N'$, we select new values $W_{N+1},\ldots
,W_{N'}$, $I_{N+1},\ldots,I_{N'}$ and $J_{N+1},\ldots,J_{N'}$ which are
independent of $W_1,\ldots,W_N$, $I_1,\ldots,I_N$ and $J_1,\ldots,J_N$,
with $W_{N+1},\ldots,W_{N'}$ independently distributed as $R(s)$, the
values of $J_{N+1},\ldots,J_{N'}$ independently and uniformly
distributed among $\{1,\ldots,m\}\setminus\{J_1,\ldots,J_N\}$ and the
values of $I_{N+1},\ldots,I_{N'}$ uniformly distributed within the
blocks $B_{J_{N+1}},\ldots,B_{J_{N'}}$, respectively. In this instance,
we take $W_1',\ldots,W_{N'}'$ and $I_1',\ldots,I_{N'}'$ to be identical
to $W_1,\ldots,W_{N'}$ and $I_1,\ldots,I_{N'}$, respectively. If
$N'<N$, then we take $(W_1',J_1'),\ldots,(W_{N'}',J_{N'}')$ to be the
(exceedence, block index) pairs that remain after randomly and
independently deleting $N-N'$ pairs from the sequence $(W_1,J_1),\ldots
,(W_{N},J_{N})$.

Let $N_0$ denote the number of exceedences of $s$ by $R_1',\ldots,R_p'$
and let $N'$ represent the number of large blocks $B_j$ in which there
is at least one exceedence of $s$ by the sequence $R_1',\ldots,R_p'$.
Then, $P(N_0\geq N')=1$. Conditional on $N_0$ and $N'$, let
$W_{N'+1}',\ldots,W_{N_0}'$ denote independent and identically random
variables, all distributed as $R(s)$, and distribute the locations
$I_{N'+1}',\ldots,I_{N_0}'$ of these exceedences independently and
uniformly over the points $\{1,\ldots,p\}\setminus\{I_1',\ldots
,I_{N'}'\}$, conditional on all of the variables $N'$, $N_0$,
$W_1',\ldots,W_{N'}'$ and $J_1',\ldots,J_{N'}'$. Take\vspace{1pt} the values of
$R_1',\ldots,R_p'$ that exceed $s$ to be the variables $W_1',\ldots
,W_{N_0}'$ and let the locations of those exceedences be the points
$I_1',\ldots,I_{N_0}'$. By construction, $W_1',\ldots,W_{N_0}'$ are
distributed as the exceedences of $s$ by $p$ independent and
identically distributed random variables distributed as $R(s)$;
conjointly, $I_1',\ldots,I_{N_0}'$ are distributed as the locations of
those exceedences and the probability that $N_0=N'=N$, $M_j\in\{0,1\}$
for each $j\in[1,m]$ and there are no exceedences of $s$ in any of the
small blocks $b_j$ for any $j\in[1,m]$ is bounded below by $1-\tau(s)$,
where $\tau(s)$ satisfies (\ref{eq2.6}); see also (\ref{eq3.20}), (\ref{eq3.21}), (\ref{eq3.23}) and~(\ref{eq3.24}).

Hence, provided that $s$ satisfies (\ref{eq3.22}), we may construct a sequence
$R_1',\ldots,R_p'$ of independent variables with the same marginal
distribution as $R_1$ and such that, with probability bounded below by
$1-\tau(s)$, the exceedences of $R_1,\ldots,R_p$ over $s$ are identical
to those of $R_1',\ldots,R_p'$. The theorem follows from this property,
(\ref{eq2.3}) and (\ref{eq3.1}), on taking $s=t$.

\section*{Acknowledgments} This paper has benefited from helpful
comments by Abba Krieger, to whom we are grateful. Research of both
authors was partially supported by an Australian Research
Council grant.

\printhistory

\end{document}